 \documentclass[11pt,a4paper]{article}
\usepackage{amsmath}
\usepackage{amsthm}
\usepackage{amssymb}
\usepackage{url}
\usepackage{amsfonts}
\usepackage{setspace}  
\usepackage{anysize}
\usepackage{enumitem}

\usepackage{amscd}        
\usepackage{graphics}
\usepackage{color}
\usepackage[all]{xy}
\usepackage{amsthm}   
\usepackage{mathrsfs}
\usepackage[all]{xy}
\usepackage{textcomp}
\usepackage{fullpage}
\usepackage{mathrsfs}
\usepackage{enumitem}
\usepackage{titlesec}
\usepackage{mathtools}

\usepackage{fancyhdr}

\titlelabel{\thetitle.\quad}

\titleformat{\subsection}[runin]{\normalfont\bfseries}{\thesubsection.}{3pt}{}
\titleformat{\subsubsection}[runin]{\normalfont\bfseries}{\thesubsubsection.}{3pt}{}

\numberwithin{equation}{section} \theoremstyle{plain}

\theoremstyle{definition}

\theoremstyle{remark}

\def\Z{\mathbb Z}
\def\N{\mathbb N}

\def\mcF{{\mathcal{F}}}
\def\C{{\mathcal{C}}}
\def\A{\mathcal{A}}

\def\rF{\mathscr{F}}
\def\rT{\mathscr{T}}

\def\rG{\mathscr{G}}

\def\ra{\rightarrow}

\def\O{{\mathcal O}}
\def\Qco{\mathfrak{Qcoh}}
\def\Mod{\mbox{-Mod}}
\def\sMod{\mbox{Mod}}

\def\fiz{\leftarrow}

\newcommand{\bfc}{\mathbf{C}}
\newcommand{\Mor}{\operatorname{Mor}}
\newcommand{\im}{\operatorname{Im}}

\newcommand{\Hom}{\mbox{Hom}}

\newcommand{\Ext}{\mbox{Ext}}
\newcommand{\id}{\operatorname{id}}

\newcommand{\Spec}{\mbox{Spec}}

\newcommand{\Ker}{\mbox{Ker}}

\newcommand{\mathcolon }{\colon\,}

\def\Qco{\mathfrak{Qco}}

\def\QcoP{\mathfrak{Qco}({\bf P^1}(R))}

\def\O{{\mathcal O}}

\def\L{{\mathscr L}}

\def\M{{\mathscr M}}

\def\Ker{\mbox{\rm Ker }}

\def\Z{\mathbb{Z}}
\def\Hom{\mbox{\rm Hom}}
\def\Ext{\mbox{\rm Ext}}
\def\Re{R[x]}
\def\Ri{R[x^{-1}]}
\def\Ro{R[x^{-1},x]}
\def\T{{\mathscr T}}
\def\Si{{S^{-1}}}
\def\B{{T^{-1}}}
\def\N{\mathbb{N}}
\def\fiz{\leftarrow}

\setenumerate{label={(\normalfont\arabic*)}} 

\begin{document}

\enlargethispage{\baselineskip}

\title{Phantom covering ideals in categories without enough projective morphisms  }

\date{}

\author{  Sergio Estrada \thanks{  Departamento de Matem\'aticas, Universidad de Murcia, 30100 Murcia, Spain. email: \texttt{sestrada@um.es}}, Pedro A. Guil Asensio\thanks{  Departamento de Matem\'aticas, Universidad de Murcia, 30100 Murcia, Spain. email: \texttt{paguil@um.es}},	
	Sinem Odaba\c{s}{\i}\thanks{  Instituto de Ciencias F\'isicas y Matem\'aticas, Universidad Austral de Chile, Valdivia-CHILE. e-mail: \texttt{sinem.odabasi@uach.cl}}
\thanks{The first and second named authors were supported by the grant MTM2016-77445-P and FEDER funds and by grant 19880/GERM/15 from the Fundaci\'on S\'eneca-Agencia
	de Ciencia y Tecnolog\'{\i}a de la Regi\'on de Murcia.  The third named author has been  supported by the research grant CONICYT/FONDECYT/Iniciaci\'on/11170394.}}

\maketitle
\renewcommand{\theenumi}{\arabic{enumi}}
\renewcommand{\labelenumi}{\emph{(\theenumi)}}

\begin{abstract}

We give sufficient conditions to ensure that the ideal $\Phi(\mathcal E)$  of $\mathcal E$-phantom maps  in a locally $\lambda$-presentable exact category $(\A, \mathcal E)$ is (special) (pre)covering ideal, where $\mathcal E$ is an exact substructure of  $(\A, \mathcal{E})$. As a byproduct, we infer the existence of various covering ideals in categories of sheaves which have a meaningful geometrical motivation. In particular we deal with a Zariski-local notion of phantom maps in categories of sheaves. We would like to point out that our approach is necessarily different from  \cite{FGHT1}, as the categories involved in most of the examples we are interested in do not have enough projective morphisms.

\end{abstract}

\vspace{0.5cm}

\section{Introduction}

Approximations by Ideals, as well as encompassing the usual Approximation Theory by Objects,  plays a key role  in certain problems such as `phantom phenomena' and `ghost phenomena' in triangulated categories, see \cite{Christensen2}, \cite{CCM08}. The authors in \cite{FGHT1} transported the `phantom phenomena' into exact category setting by means of additive subfunctors $\mathcal F$ of $\Ext$ in $(\A,\mathcal{E})$, see Definition \ref{relative phantom}. For a given such functor $\mathcal{F}$, the ideal of $\mathcal{F}$-phantoms is denoted by $\Phi(\mathcal{F})$. The main contribution of the aforementioned paper is that if $(\A,\mathcal{E})$ has  enough injective objects and projective morphisms (see  (\ref{projective morphism})), then  an ideal $\mathcal I$ of $\A$ is special precovering (see Definition \ref{idealcover}) if and only if $\mathcal I$ is of the form $\Phi(\mathcal F)$, for some additive subfunctor $\mathcal F\subseteq \Ext$ with enough injective morphisms.

However, there are many natural exact (even abelian) categories which have no  enough projective morphisms (see Proposition \ref{nohay}),  so the previous theorem can no longer be applied to infer that an ideal of the form $\Phi(\mathcal F)$ is  special (pre)covering.

One of the main purposes of this paper is to develop new methods, necessarily different from those in \cite{FGHT1}, to remove the hypothesis of the existence of enough projective morphisms in \cite[Theorem 17]{FGHT1} to assure that ideals of the form  $\Phi(\mathcal F)$ for certain  additive subfunctors $\mathcal{F}$ are special (pre)covering. Regarding this, the main result of Section 3   is the following:

\medskip\par\noindent
{\bf Theorem \ref{main} - Corollary \ref{spec}.} Let $(\A,\mathcal{E})$ be an exact category with $\A$ locally presentable.  If  $\mathcal E'$ is an exact substructure  of $(\A,\mathcal{E})$ which   is closed under direct limits, then $\Phi(\mathcal{E}')$ is a covering ideal. In addition, if $\A$ has   enough $\mathcal{E}'$-phantom morphisms, and $\mathcal{E}'$ has enough injectives,  then $\Phi(\mathcal{E}')$ is a special covering ideal.

\medskip\par\noindent

In Section 4, we focus on the category $\Qco(X)$ of quasi-coherent sheaves in which  \cite[Theorem 17]{FGHT1} can not be directly applied unless $X$ is an affine scheme.

\medskip\par\noindent

{\bf Proposition \ref{nohay}.} There are no nonzero  projective morphism in $\Qco({\bf P^1}(R))$.

\medskip\par\noindent

One immediate consequence of this result is the following.

\medskip\par\noindent

 {\bf Theorem. \ref{nohay.phantom}.}
There are no non-trivial phantom maps in $\Qco({\bf P^1}(R))$.

\medskip\par\noindent

On the other hand, for any scheme $X$,  the category $\Qco(X)$ is known to be Grothendieck with a closed monoidal structure $\otimes$. This fact together with Definition \ref{stalkwise} leads to the following inclusion of exact structures  in $\Qco(X)$
$$\mathcal{E}_{st} \subseteq \mathcal{E}_{\otimes}.$$
 From \cite{EEO} and $\cite{EGO15}$, we know that they have enough injectives. So we have the following result:

 \medskip\par\noindent

 {\bf Corollary \ref{covering.ideals.qco}.} If $X$ is a quasi-compact and semi-separated scheme, then $\Phi(\mathcal{E}_{st})= \Phi(\mathcal{E}_{\otimes})$ is a special covering ideal which has a Zariski-local property.

 \medskip\par\noindent

\section{Preliminaries}

In this section, we recall certain concepts and facts which are pertinent to our study.

\subsection{Exact categories}\cite{Quillen}. Let $\A$ be an additive category, and $\mathcal{E}$ be a class of  kernel-cokernel pairs $(i,p)$ in $\A$
$$\xymatrix{A\ar[r]^i & B \ar[r]^p &C},$$
that is,  $i$ is the kernel of $p$, and $p$ is the cokernel of $i$. The morphisms $i$ and $p$ are said to be  an \emph{admissible monic} and an \emph{admissible epic} in $\mathcal{E}$, respectively.  The class $\mathcal{E}$ is said to form an \textit{exact structure on $\A$} if it is closed under isomorphisms and satisfies the following:
\begin{itemize}
	\item[E0)] For every object $A$ in  $\A$, the identity morphism $\id_A$ is both an admissible monic and an admissible epic in $\mathcal{E}$.
	\item[E1)] The classes of admissible monics and admissible epics in $\mathcal{E}$ are closed under compositions.
	
	\item[E2)] The pushout and pullback of a kernel-cokernel pair in $\mathcal{E}$ along an arbitrary morphism exist and they are a kernel-cokernel pair in $\mathcal{E}$.
		\end{itemize}

An \textit{exact category} $(\A, \mathcal{E} )$ consists of an additive category $\A$, together with an  exact structure $\mathcal{E}$ on $\A$. Elements of $\mathcal{E}$ are called \textit{short exact sequences}. From now on, by such arrows
$\xymatrix{\ar@{^{(}->}[r]&}$ and  $\xymatrix{\ar@{->>}[r]&}$ we mean to  admissible monics and admissible epics in $\mathcal{E}$, respectively.

For any additive category $\A$,  we let $\mathcal{E}_{sp}$ denote the smallest exact structure of splitting short exact sequences in $\A$.

\vspace{3mm}

The following lemma plays a  crucial role in the proof of our main result, Theorem \ref{main}.

\subsection{Lemma}\label{Obscure}\cite[Corollary 7.5, Proposition 7.6]{Buhler}.
Let $(\A,\mathcal{E})$ be an exact category. The following are equivalent:
\begin{itemize}
\item[i)] The additive category $\A$ is weakly idempotent complete.

\item[ii)] Every retraction is an admissible epic in $\mathcal{E}$.
\item[iii)] If the composition $g \circ f$ is an admissible
epic in $\mathcal{E}$, then $g$ is an admissible epic in $\mathcal{E}$.

\end{itemize}

\subsection{Projectively (injectively) generated exact structures.}
Let $(\A,\mathcal{E})$ be an exact category, and $\mathcal{S}$ be a class of objects in $\A$. The exact structure $\mathcal{E}$ is said to be \textit{projectively (injectively) generated by the class} $\mathcal{S}$ if $\mathcal{E}$ is the largest exact structure on $\A$ for which each $S\in \mathcal{S}$ is $\mathcal{E}$-projective ($\mathcal{E}$-injective) (for definitions of $\mathcal{E}$-projective and $\mathcal{E}$-injective objects, see \cite[Section 11]{Buhler}).
\subsection{Remark.} If $(\A,\mathcal{E})$ has enough  projective (injective) objects, then the exact structure $\mathcal{E}$  is  projectively (injectively) generated by the class of projective (injective) objects, see \cite[Exercise 11.10]{Buhler}. However, the converse doesn't hold in general: Let $\A$ be an  abelian category  without enough projectives, and $\mathcal{E}$ be  the exact structure of usual short exact sequences. It is projectively generated by the set $\{0\}$, however, has no  enough projectives.

\subsection{} For a fixed exact category $(\A,\mathcal{E})$ and for given $A,C\in \A$, we let  $\Ext(C,A)$ denote  the class of  all  representatives of isomorphism classes of short exact sequences in $\mathcal{E}$ beginning with $A$ and ending with $C$.

An \textit{exact substructure}  $\mathcal{E}'$ of $\A$ is just an exact structure on $\A$ with $\mathcal{E}' \subseteq \mathcal{E}$. In such case, we shall denote by $\Ext_{\mathcal{E}'}(C,A)$ the group of all representatives of isomorphism classes of short exact sequences in $\mathcal{E}'$.  Note that any exact substructure  $\mathcal{E}'$ on $\A$ yields the additive subfunctor $\Ext_{\mathcal{E}'}$ of $\Ext$. However, the converse isn't true in general, see \cite[pg. 651]{DRSSK99}.

\subsection{Presentable categories.}\label{present}\cite{AR}. Let  $\lambda$ be a  regular cardinal. An object $A$ in a category $\A$ is said to be \emph{$\lambda$-presentable} if the functor $\Hom_{\A}(A,-)$
preserves $\lambda$-directed colimits. The category $\A$ is called
\emph{ locally $\lambda$-presentable} if it is cocomplete, and there is a
set $\mathcal S$ of $\lambda$-presentable objects in $\A$ such that any other object in $\A$ is a $\lambda$-directed colimit of objects in
$\mathcal S$.
For short, $\aleph_0$-directed colimits will be just called \emph{direct
	limits;} locally $\aleph_0$-presentable categories,  \emph{locally finitely presented}; and $\aleph_0$-presentable objects, \emph{finitely presented}.

Two of the main properties of locally presentable categories which will be pertinent in the proof of Proposition \ref{cover} are  that they are well-powered (see \cite[Remark 1.56]{AR}) and every object is presentable (see \cite[Proposition 1.16]{AR}).

\subsection{Pure quotients.}\label{purequotient}\cite{AR04}. Let $\lambda $ be a regular cardinal. A morphism $f\!:\! A \rightarrow B$ in any category $\A$ is called \emph{$\lambda$-pure quotient} if for any morphism $v:F \rightarrow B$
with $F$  $\lambda$-presentable, there is a morphism $g: F \rightarrow A$ such that $v=f \circ g$

$$\xymatrix{& F \ar[d]^v \ar@{-->}[ld]_g\\
	A\ar[r]^f & B.}$$

If $\A$ is $\lambda$-accessible with pullbacks, we infer from  Ad\'{a}mek and Rosick\'{y} \cite[Proposition 3]{AR04} that a morphism in $\A$ is $\lambda$-pure quotient if and only if it is a $\lambda$-directed colimit of retractions.
In particular,  if $\A$ is locally $\lambda$-presentable, then it is complete (see \cite[Remark 1.56]{AR}), and therefore,  has pullbacks.

\subsection{Remark.} If $\A$ is an additive and  finitely accessible category, we already know from \cite{CB} that there is the so-called \textit{pure exact structure}, denoted by $\mathcal{E}_{\aleph_0}$, which has enough projectives, and is closed under direct limits. An admissible epic in $\mathcal{E}_{\aleph_0}$ is necessarily an $\aleph_0$-pure quotient.

On the other hand, if $\A$ is abelian and $\lambda$-accessible for some regular cardinal $\lambda$,     it is not difficult to show that   $\lambda$-pure quotients give rise to the exact structure $\mathcal{E}_{\lambda}$ which is also  projectively generated by $\lambda$-presentable objects.

\subsection{Morphism category.} Given any category $\A$, we let   $\Mor(\A)$ denote the category of morphisms of $\A$. It is well-known that the category $\Mor(\A)$ inherits many categorical properties
from $\A$: if $\A$ is additive, so is $\Mor(\A)$; if $\A$ is  $\lambda$-accessible, so is $\Mor(\C)$; if $\mathcal{E}$ is an exact structure on $\A$, then the class $\Mor(\mathcal{E})$ of morphisms between short exact sequences in $\mathcal{E}$ leads to an exact structure on $\Mor(\A)$.

\subsection{Monoidal categories and purity.} A monoidal category is a category $\A$ equipped with a bifunctor
$$-\otimes - \mathcolon  \quad \A \times \A \rightarrow \A$$
and a unit object $I$, subject to certain coherence conditions which ensure that all relevant diagrams commute. The monoidal category $(\A,\otimes, I)$ is said to be  \textit{symmetric} if,  for every pair of objects $A, B$  in $\A$, there is an isomorphism $A \otimes B \cong B \otimes A$, which is natural on both $A$ and $B$. The monoidal structure $\otimes$ on $\A$ is said to be \emph{closed} if for each $A\in \A$, the functor $- \otimes A\mathcolon  \A \ra \A$ has a right adjoint $[A,-]\mathcolon  \A \ra \A$. For the whole axioms  and examples of (closed) symmetric monoidal categories, see \cite{Kelly}.

Let  $f\mathcolon X \rightarrow Y$ be a monomorphism  in a closed symmetric  monoidal category $(\A,\otimes,I)$.  The author in  \cite{Fox} calls it   \emph{$\otimes$-pure} if for every $Z \in \C$, $f\otimes Z$ remains monic. Following  \cite{EGO15}, we will call it  \emph{geometrical pure}. If $\A$ is abelian, the class $\mathcal{E}_{\otimes}$ of geometrical pure short exact sequences forms an exact structure on $\A$. In \cite{EGO15}, it was proved that if, in addition, $\A$ is Grothendieck, then $(\A, \mathcal{E}_{\otimes})$ is an exact category of Grothendieck type (\cite[Section 3]{Sto12}),  and therefore,  it has enough injectives.

\subsection{Remark.}The geometrical pure exact structure recovers  many of the known exact structures. For instance, the usual purity in $R\Mod$; the componentwise purity, as well as,  the categorical purity in the category $\bfc(R)$ of complexes of $R$-modules; the usual purity on stalks in the category $\O_X \Mod$ of $\O_X$-modules; or the stalkwise purity in the category $\Qco(X)$ of quasi-coherent sheaves over a quasi-separated scheme $X$.

\section{Covering ideals of relative phantom morphisms}

 \subsection{Ideals.} An \emph{ideal} $\mathcal{I}$ of a preadditive category $\A$ is just an additive subfunctor of $\Hom$, equivalently,  a family  $\{\mathcal{I}(A,A')\}_{A,A' \in \A}$ of subgrops of morphisms in $\A$ which  are closed under left and right compositions by other morphisms.

 As an easy observation,   a morphism  $f\mathcolon  A\rightarrow A'$ in $\A$ belongs to $\mathcal{I}$ if and only if   the natural transformation $$\Hom(-,f)\mathcolon  \Hom_{\A}(-,A)\rightarrow \Hom_{\A}(-,A')$$ factors through $\mathcal{I}(-,A')$.

 \subsection{(Pre)covering and (pre)enveloping ideals.}\label{idealcover}Let $\mathcal{I}$ be a class of morphisms  in a  category   $\A$. An \emph{$\mathcal{I}$-precover} of an object $X\in \A$ is  a morphism $f\mathcolon  A\rightarrow X$ in $\mathcal{I}$ such that any other morphism $f'\mathcolon  A'\rightarrow X\in \mathcal{I}$  factors through $f$.  It is said to be an $\mathcal{I}$ \emph{cover} if every morphism $g$ with $f \circ g=f$ is an isomorphism.

 In addition, if $\A$ is endowed with an exact structure $\mathcal{E}$, then an  $\mathcal{I}$-precover $f\mathcolon   A \rightarrow X$  of $X$ is said to be \emph{special} if $f$ is an admissible epimorphism in $\mathcal{E}$ and it fits in a pushout diagram
 $$
 \xymatrix{
 K\ar@{^{(}->}[r]\ar[d]_g& A'\ar@{->>}[r]\ar[d]  &X\ar@{=}[d]\\
E\ar@{^{(}->}[r] &  A\ar@{->>}[r]^f &   X }
 $$
 with  $\Ext(i,g)=0$ for every $i \in \mathcal{I}$. The class $\mathcal{I}$ is said to be   \emph{(special) (pre)covering} if every object in $\A$ has a (special) $\mathcal{I}$-(pre)cover.

  As an aside,   one can easily observe that if $\mathcal{I}$ is an ideal in a preadditive category $\A$, then being $f$ an $\mathcal{I}$-precover is equivalent to   $\im \Hom(-,f)=\mathcal{I}(-,X)$.

The notions of (special) $\mathcal{I}$-(pre)envelope are defined dually.

  \subsection{Remark.}\label{object-ideal}
  In Relative Homological Algebra, approximations are usually constructed from a (pre)covering and/or (pre)enveloping   class of objects. On the other hand, Approximation Theory by morphisms/ideals naturally generalizes  the usual Approximation Theory by objects. Indeed, for   a given   class $\mcF$ of objects in   $\A$, $\mcF$ is (pre)covering ((pre)enveloping)  if and only if the class $\Hom(\mcF,\textrm{-})$ $(\Hom(\textrm{-}, \mcF))$ of all morphisms with domain (codomain) in $\mcF$ is (pre)covering ((pre)enveloping).

  If, in addition,  $\A$ is  an additive category,
  one can easily show that an additive subcategory $\mathcal{F}$ of $\A$ is precovering (preenveloping)  if and only if
  the \textit{object ideal } $\mathcal{I}(\mcF)$ \textit{generated by} $\mcF$ (see \cite[Proposition 1]{FH16})
 is a precovering (preenveloping) ideal.

 \subsection{Relative phantom morphisms}\label{relative phantom}\cite[Definition 6]{FGHT1}. Let $(\A,\mathcal{E})$ be an exact category. Given an exact substructure $\mathcal{E}'$, a morphism $f\mathcolon  A \rightarrow B$ is called  \emph{$\mathcal{E}'$-phantom} if $\Ext(f,\textrm{-})$ transforms short exact sequences in $\mathcal{E}$ into  short exact sequences in $\mathcal{E}'$, that is, $$\Ext(f,\textrm{-})\mathcolon  \Ext(B,\textrm{-})\rightarrow \Ext_{\mathcal{E}'}(A, \textrm{-}).$$ The class  of $\mathcal{E}'$-phantom morphisms, denoted by $\Phi(\mathcal{E}')$,  constitutes  an ideal in $\A$.  We say that $\A$ \textit{has enough $\mathcal{E}'$-phantom morphisms} if for every object $B\in \A$ there exists an $\mathcal{E}'$-phantom morphism $f:A \rightarrow B$ which is  an  admissible epic in $\mathcal{E}$.

  \subsection{Example.}\label{projective morphism} Let $f$ be a morphism in  an exact  category $(\A, \mathcal{E})$. It is called    \emph{projective} if $f \in \Phi(\mathcal{E}_{sp})$. Note that if $(\A,\mathcal{E})$ has enough projective objects, then  $f$ is projective if and only if it factors through a projective object. $\A$ said to have \textit{have enough projective morphism} if it has enough  $\mathcal{E}_{sp}$-phantom morphisms.

 \subsection{Remark.}
 In certain cases, in order to check whether or not a morphism $f\mathcolon A \rightarrow B$ in an exact category is $\mathcal{E}'$-phantom,  one need not take the pullback along \emph{all} short exact sequences: for instance, suppose that $(\A,\mathcal{E})$ is an efficient exact category with a generator (see \cite{Sto12}). By \cite[Proposition 5.3]{Sto12}, for every object $B \in \A$ there is   a set of short exact sequences
 $$\mathbb{E}_i: \quad \xymatrix{ K_i \ar@{^{(}->}[r] & T_i \ar@{->>}[r] & B}$$
 in $\mathcal{E}$  indexed by $i\in I_B$, such that every short exact sequence $\mathbb{E}$ in $\mathcal{E}$ ending with  $B$ is a pushout of   some of them. That is, $\mathbb{E} \cong g \mathbb{E}_i$, for some $i\in I_B$ and some morphism $g$. Therefore, 
 for a given  exact substructure  $\mathcal{E}'$,  $f\mathcolon A\rightarrow B $ in $\A$ is $\mathcal{E'}$-phantom if and only if $\mathbb{E}_i f$ belongs to $\mathcal{E}'$, for every $i\in I_B$.

 In addition, if $\A$ has enough projective morphisms, one can work with  just one  short exact sequence
$$\mathbb{E}: \quad \xymatrix{ K \ar@{^{(}->}[r] & P \ar@{->>}[r] & B}$$
 whose admissible epic is a projective morphism.

  \subsection{Remark.}Inspired from phantom morphisms in the stable module category $kG\mbox{-}\underline{\sMod}$ of modules over $kG$, where $k$ is a field and $G$ is a finite group (for example, see  \cite{Gnacadja}) the author in  \cite{Herzog} defines phantom morphisms in the module category $R\Mod$  over an arbitrary ring $R$. He also proves that a  morphism $f$ is phantom in $R \Mod$  if and only if $f$ is $\mathcal{E}_{\aleph_0}$-phantom,  where $\mathcal{E}_{\aleph_0}$ is the proper class of all pure short exact sequences in $R\Mod$. This observation had two consequences: on one side, it leaded to the authors in \cite{FGHT1} to define `relative' phantom morphisms in any exact category; on the other side, it showed in this abelian framework the close relation between phantom and purity, which was already revealed in the triangulated setting in \cite{Christensen1} and  \cite{Christensen2}.

 \subsection{Terminology.}
 Following the previous motivation, in case $\A$ is locally finitely presented, then we call a morphism in  $\Phi(\mathcal{E}_{\aleph_0})$  just \textit{phantom} without the prefix.

  \vspace{5mm}

Before going further, we prove the following:

\subsection{Proposition.}\label{cover} Let $\mathcal{I}$ be a class (not necessarily ideal) of morphisms in an additive locally presentable category $\A$ which is closed under direct limits. If an object $X$ of $\A$ has an $\mathcal{I}$-precover, then it has an $\mathcal{I}$-cover.

\vspace{5mm}

The previous result is known to be true  for a class of objects in a Grothendieck category  (\cite[Theorem 1.2]{ElBashir}). In Proposition \ref{cover}, we extend it to any class of morphisms in a locally presentable category. We basically adopt the proof given in \cite{Eno81} for $R \Mod$ using techniques of subobjects (\cite[Section IV.1]{Ste}) and presentability of objects as mentioned in (\ref{present}). We aim to prove Proposition \ref{cover} in the next three lemmas. We assume that  $\mathcal{I}$ is  a class of morphisms in a locally presentable category  $\A$ which is closed under direct limits.

\subsection{Lemma.}\label{lemma1}Suppose that an object $X$ in $\A$ has an $\mathcal{I}$-precover. Given any morphism $f:F \rightarrow X$ in $\mathcal{I}$, there exists a factorization $f=t \circ s$ satisfying:
\begin{itemize}
	
	\item[a)]$t$ is an $\mathcal{I}$-precover of $X$.
	\item[b)] If $t'$ is any other $\mathcal{I}$-precover of $X$, then for  any factorization $t =t' \circ h$, $\Ker(s)$ and $\Ker(h\circ s)$ together with the canonical inclusions are equivalent subobjects of $F$. In particular, $\Ker(s)\cong \Ker(h \circ s)$.
	
\end{itemize}

 \begin{proof}
Assume, for the sake of contradiction,  that there is a morphism $f\mathcolon F \rightarrow X$ in $\mathcal{I}$ which doesn't satisfy the conditions given. From \cite[Remark 1.56]{AR}, $\A$ is well-powered. So the class  $\mbox{Subob}(F)$ of subobjects of $F$ is in fact a set. We choose an infinite ordinal $\gamma$ whose cardinality is strictly bigger than than the cardinality of     $\mbox{Subob}(F)$.  
 By transfinite induction, we construct a $\gamma$-directed continuous system $\{A_{\alpha}:\ s_{\alpha \beta}\}_{\alpha< \beta \leq \gamma}$ together with a cone morphism $\{t_{\alpha}\mathcolon A_{\alpha} \rightarrow X\}_{\alpha \leq \gamma}$ satisfying:
 \begin{itemize}
 	\item[i)] $t_0=f$.
 	\item[ii)] For every ordinal $\alpha >0$, $t_{\alpha}$ is an $\mathcal{I}$-precover of $X$.
 	\item[iii)] For every ordinals $0 < \alpha < \beta \leq \gamma$, $\Ker (s_{0 \alpha})  \subsetneqq  \Ker{s_{0 \beta}}$.
 	
 \end{itemize}
 	 Firstly, let $t_0:=f$.  Since $X$ has an $\mathcal{I}$-precover, say $t_1\mathcolon A \rightarrow X$, $t_0$ has a factorization $t_0= t_1 \circ s_{01}$. Now, let $0< \beta< \gamma$ be  an ordinal such that the family $\{A_{\alpha}:\ s_{\alpha' \alpha}\}_{\alpha'<\alpha \leq \beta}$ together with a cone morphism $\{t_{\alpha}\mathcolon A_{\alpha} \rightarrow X\}_{\alpha \leq \beta}$ is constructed. By assumption,  there exists an $\mathcal{I}$-precover $t_{\beta+1}$ of $X$ and a
 	 factorization $t_{\beta}=t_{\beta+1} \circ s_{\beta \beta+1}$ such that
 	 $\Ker (s_{0 \beta})  \subsetneqq  \Ker(s_{\beta \beta+1}\circ s_{0 \beta} )$. For any ordinal $\alpha \leq \beta $, we let $s_{\alpha \beta+1}:= s_{\beta \beta+1} \circ s_{\alpha \beta}$, and we have $\Ker(s_{0\alpha}) \subsetneqq \Ker (s_{0\beta})  \subsetneqq  \Ker(s_{0\beta+1} ) $.

 	 If $\beta \leq \gamma$ is a limit ordinal, then we simply let $t_{\beta}:= \varinjlim_{\alpha < \beta} t_{\alpha}$. By assumption, it is a morphism in $\mathcal{I}$, and one can easily show that in fact, it is  an $\mathcal{I}$-precover of $X$. For any ordinal $\alpha < \beta$, we let $s_{\alpha \beta}:\ A_\alpha \rightarrow A_\beta=\varinjlim_{\alpha < \beta} A_\alpha $ be the canonical morphism. If  for some ordinals $0<\alpha < \beta$, the subobjects  $\Ker (s_{0 \alpha}) \subseteq \Ker({s_{0 \beta}})$ of $F$ were  equivalent, that is  $  \Ker({s_{0 \beta}}) \subseteq \Ker (s_{0 \alpha})$, then
 	$$  \Ker (s_{0 \alpha+1}) \subseteq  \Ker({s_{0 \beta}}) \subseteq \Ker (s_{0 \alpha})$$
 	and therefore, $ \Ker (s_{0 \alpha}) \subseteq  \Ker (s_{0 \alpha+1})$ would be  equivalent subobjects of $F$, which leads to a contradiction. So
 	$  \Ker({s_{0 \alpha}}) \subsetneqq \Ker (s_{0 \beta})$.

 	As a result of the construction, there exists  a  chain of distinct subobjects of $F$
 	$$\Ker (s_{0 1}) \subsetneqq \Ker({s_{0 2}})   \subsetneqq \cdots \subsetneqq \Ker (\varphi_{0 \gamma}).$$
 But, the cardinality of $\mbox{Subob}(F)$ is strictly less than $ \gamma$,   which leads to  a contradiction.

\end{proof}

 \subsection{Lemma.}I\label{lemma2}f an object $X$ of $\A$ has an $\mathcal{I}$-precover, then it has an $\mathcal{I}$-precover  $\varphi$ such that for any factorization $\varphi= h \circ f$ where $h$ is an $\mathcal{I}$-precover of $X$, $f$ is a monomorphism.

 \begin{proof}
 	Let $\varphi_0: F \rightarrow X$ be an $\mathcal{I}$-precover of $X$. Using Lemma \ref{lemma1}, one can construct a directed system $\{F_n: \varphi_{n m}\}_{0 \leq n<m \in \N}$ with a cocone $\varphi_n: F_n \rightarrow X$ satisfying that  for every natural number $0 \leq n $, $\varphi_n$ is an $\mathcal{I}$-precover of $X$, and  the factorization $\varphi_n =\varphi_{n+1} \circ \varphi_{nn+1}$ satisfies Lemma \ref{lemma1}.

 	Let $\varphi:= \varinjlim_{0 \leq n} \varphi_n$.  We claim that the morphism $\varphi$ is the desired one. Clearly, it is an $\mathcal{I}$-precover of $X$. Suppose that there exists a a factorization $\varphi= h \circ f$ where $h$ is an $\mathcal{I}$-precover of $X$. In order to show that $f$ is a monomorphism, it sufficies to show $\Ker(f)=0$. Let $k\mathcolon K \rightarrow \varinjlim_{0 \leq n} F_n$ be the  kernel of $\varphi$.  Note that $K$ is $\lambda$-presentable for some regular cardinal $\lambda$,  see \cite[Proposition 1.16]{AR}. So there exists a morphism $k'\mathcolon K \rightarrow F_n$  for some $n \geq 0$ such that  $k=g_n \circ k'$, where $g_n: F_n \rightarrow \varinjlim_n F_n$ is the canonical morphism. But $g_n=g_{n+1} \circ \varphi_{nn+1}$ and
 	$$0=f \circ k=f \circ g_{n+1} \circ \varphi_{nn+1} \circ k'$$
 	So the morphism $k'$ is factorized over $\Ker (f \circ g_{n+1} \circ \varphi_{nn+1})$. By construction, the factorization $\varphi_n =\varphi_{n+1} \circ \varphi_{nn+1}$ satisfies Lemma \ref{lemma1}, and therefore, $\Ker(\varphi_n) \cong \Ker (f \circ g_{n+1} \circ \varphi_{nn+1})$. This implies that $\varphi_{nn+1} \circ k'=0$,  and hence, $k=0$.
 	
 \end{proof}

 \subsection{Lemma.}\label{lemma3}If $\varphi: A \rightarrow X$ is an $\mathcal{I}$-precover satisfying Lemma \ref{lemma2}, then it is an $\mathcal{I}$-cover.
 \begin{proof}
 	Suppose that $\varphi=\varphi \circ f$ for some morphism $f\mathcolon A\rightarrow A$. By assumption, it is a monomorphism, and therefore, $(A,f)$ is a subobject  of $A$.
 	
  Let   $\gamma$  be an infinite ordinal whose cardinality  is strictly bigger than the cardinality of $\mbox{Subob}(A)$. By transfinite induction, we construct a directed system $\{A_{\alpha}, f_{\alpha \beta}\}_{\alpha < \beta \leq \gamma}$, satisfying:
 	\begin{itemize}
 		\item[i)] For every ordinal $\alpha \leq \gamma$, $A_{\alpha}=A$.
 		\item[ii)] For every ordinal $\alpha <\gamma$, $f_{\alpha \alpha+1}=f$,
 		\item[iii)] For every ordinals $\alpha < \beta$,  $\varphi=\varphi \circ f_{\alpha \beta}$.
 	\end{itemize}
 	
 	Firstly, let $A_0:=A$.  Suppose that such family is constructed for every ordinal $\leq \beta$, for some ordinal $\beta < \gamma$. Then  	we let  $A_{\beta+1}:=A$, $f_{\beta\beta +1}:=f$, and $f_{\alpha \beta+1}:= f \circ f_{\alpha\beta}$,  for every ordinal $\alpha \leq \beta$.

 	If $\beta \leq \lambda$ is a limit ordinal, consider the morphism $\varinjlim_{\alpha< \beta} \varphi\mathcolon \varinjlim_{\alpha < \beta } A_{\alpha} \rightarrow X$. By assumption, it is a morphism in $\mathcal{I}$. So using the fact that $\varphi$ is an $\mathcal{I}$-precover, there exists a factorization  $\varinjlim_{\alpha< \beta} \varphi= \varphi \circ s$.
 	Then for every ordinal $\alpha < \beta$, we let $f_{\alpha \beta}$ denote the composition of $s$ with the canonical morphism $A_{\alpha} \rightarrow \varinjlim_{\alpha < \beta } A_{\alpha}.$ As a result, we get the desired system.
 	
 	On the other side, by assumption, for every $\alpha < \beta \leq \gamma$,  the morphism $f_{\alpha \gamma}:A\rightarrow A$ is a monomorphism  and  	$f_{\alpha \gamma}=f_{\beta \gamma} \circ f_{\alpha \beta}$. It gives rise to  a $\gamma$-length chain of subobjects of $A$
 	$$(A,f_{0\gamma}) \subseteq (A,f_{1\gamma}) \subseteq \cdots \subseteq (A, f_{\gamma \gamma}).$$
 But  $A$ has less than $\mid \gamma \mid$ subobjects, so there must be at least two different ordinals $\alpha < \beta \leq \gamma$ such that the subobjects $(A, f_{\alpha \gamma}) \subseteq (A, f_{\beta \gamma})$ of $A$ are equivalent, that is,  $(A,f_{\beta \gamma}) \subseteq (A ,f_{\alpha\gamma})$. This implies that
 	$$(A,f_{\alpha+1 \gamma}) \subseteq (A,f_{\beta \gamma}) \subseteq (A ,f_{\alpha\gamma}),$$
 	and therefore, the subobjects $(A ,f_{\alpha\gamma}) \subseteq  (A ,f_{\alpha+1\gamma})   $ of $A$ are equivalent. 
Using property of monomorphism and $f_{\alpha+1 \gamma} \circ f=f_{\alpha \gamma}$, we conclude that  $f$  is an isomorphism.
 \end{proof}

\subsection{Lemma.}\label{krau}\cite[Theorem 4]{Krause}.
 	Let $\A$ be an additive and locally presentable category, and  let $\mathcal B$ be an additive subcategory of $\mathcal A$ which is closed under direct limits. If $\mathcal B$ is closed under $\lambda$-pure subobjects or $\lambda$-quotients for some regular cardinal $\lambda$, then it is a precovering class in $\mathcal A$.

\subsection{Remark.}Combining Remark \ref{object-ideal} with Proposition \ref{cover}, we can say now that the class $\mathcal{B}$ in Lemma \ref{krau} is in fact a covering class.

 \vspace{3mm}

From now on,  $\A$ stands  for an exact category with a fixed exact structure $\mathcal{E}$.
\subsection{Theorem.}\label{main}
Suppose that $\A$ is locally $\lambda$-presentable for some regular cardinal $\lambda$.   If  $\mathcal E'$ is an exact substructure  of $\A$ which   is closed under direct limits, then $\Phi(\mathcal{E}')$ is a covering ideal.

\begin{proof}
Firstly, note that 	since $\A$ is additive and  locally presentable category, it is weakly idempotent complete. Therefore, by Lemma \ref{Obscure}, any retraction is an admissible epic in $\mathcal{E}'$.   Using the fact given in (\ref{purequotient}) and, the assumption on $\mathcal{E'}$,  we deduce that any $\lambda$-pure quotient is an admissible epic in $\mathcal{E}'$.
	
	On the other hand,  the ideal $\Phi(\mathcal{E}')$ is an additive subcategory of $\Mor(\A)$. Due to the fact that a pullback diagram is a finite limit and finite limits commutes with colimits, one can easily deduce  that the  ideal $\Phi(\mathcal{E}')$ is closed under direct limits.  From Proposition \ref{cover}, it suffices to show that $\Phi(\mathcal{E}')$ is a precovering ideal in $\A$. We claim  that $\Phi(\mathcal{E}')$ is closed under $\lambda$-pure quotients in $\Mor(\A)$. For that, let $f\mathcolon A\rightarrow B$ be a morphism in $\Phi(\mathcal{E}')$ and $a\mathcolon  f \rightarrow f'$ be  a $\lambda$-pure quotient in $\Mor(\A)$
	$$\xymatrix{A\ar[r]^{a_0} \ar[d]_f & A' \ar[d]^{f'} \\
		B\ar[r]_{a_1} & B' .
	}
	$$
	Since $\lambda$-presentable objects in $\Mor(\A)$ are precisely morphisms between  $\lambda$-presentable objects,   $a_0$ and $a_1$ are $\lambda$-pure quotients in $\A$.
	Besides, by assumption we have
	$$\Ext(a_1 \circ f, \textrm{-})= \Ext(f' \circ a_0, \textrm{-})\mathcolon \Ext(B', \textrm{-}) \rightarrow \Ext_{\mathcal{E}'}(A, \textrm{-}).$$
	Then, for any short exact sequence $\mathbb{E} \in \mathcal{E}$, $(\mathbb{E}f')a_0=\mathbb{E}(f' \circ a_0)$ belongs to  $\mathcal{E}'$,
	$$\xymatrix{
		\mathbb{E}(f'\circ a_0): & X \ar@{^{(}->}[r]\ar@{=}[d]&Y''\ar[d]\ar@{->>}[r]^g& A\ar[d]^{a_0}\\
		\mathbb{E}f': & X \ar@{^{(}->}[r]\ar@{=}[d]&Y'\ar[d]\ar@{->>}[r]& A'\ar[d]^{f'}\\
		\mathbb{E}: & X \ar@{^{(}->}[r]&Y\ar@{->>}[r]& B'.}
	$$
	As $a_0$ and $g$ are admissible epics in  $\mathcal{E}'$,  so is  $a_0 \circ g$. By Lemma \ref{Obscure},  $\mathbb{E}f' \in \mathcal{E}'$. Therefore, we infer from Lemma \ref{krau} that $\Phi(\mathcal{E}')$ is a precovering class in $\Mor(\C)$, from which it follows easily that it is a precovering ideal in $\A$.
	
\end{proof}

  \subsection{}\label{push1}Let $\xymatrix{K \ar[r]^k & B\ar[r]^{\varphi} & A}$ be a sequence, that is, $\varphi \circ k=0$.  Let
 $\xymatrix{K \ar@{^{(}->}[r]^e & E \ar@{->>}[r]^a & K'}$ be a short exact sequence in $\mathcal{E}'$. Then there exists the following pushout diagram
 \begin{equation}\label{d31}
 \xymatrix{
 	K \ar@{^{(}->}[r]^e\ar[d]_k&E \ar@{->>}[r]^a\ar[d]^{k'}&K'\ar@{=}[d]\\
 	B\ar@{^{(}->}[r]^{v}\ar[d]_{\varphi}& B'\ar@{->>}[r]^{a'}\ar[dl]^-{\varphi'}& K'\\
 	A& & }
 \end{equation}
 with $\varphi' \circ k'=0$.
 \subsection{Lemma.}\label{pushoutp}
 Consider the pushout diagram given in (\ref{d31}). If $\varphi$ is an $\mathcal{E}'$-phantom morphism. Then $\varphi'$ is $\mathcal{E}'$-phantom, as well.

 \begin{proof}

 	Let
 	$$\mathbb{E}:\ \xymatrix{Y \ar@{^{(}->}[r] & X \ar@{->>}[r] &A}$$
 	be a short  exact sequence in $\mathcal{E}$. Then there exists a commutative diagram
 	$$\xymatrix@=.5cm{
 		&&&&K\ar@{=}[rr]\ar[dd]^<<<<t\ar@{^{(}->}[dl]^e&&K\ar[dd]^k\ar@{^{(}->}[dl]^e\\
 		&&&E\ar@{=}[rr]\ar[dd]^<<<<{t'}&&E\ar[dd]^<<<<{k'}&\\
 		(\mathbb{E}\varphi):&&Y\ar@{^{(}->}[rr]^<<<<<{y''}\ar@{=}[dl]\ar@{=}[dd]&&X''\ar@{->>}[rr]\ar[dl]^{v'}\ar[dd]&&B\ar[dd]^-{\varphi}\ar@{^{(}->}[dl]^{v}\\
 		(\mathbb{E}\varphi'):&Y\ar@{^{(}->}[rr]^>>>>>{y'}\ar@{=}[dd]&&X'\ar@{->>}[rr]\ar[dd]&&B'\ar[dd]^<<<<{\varphi'}&\\
 		&&Y\ar@{^{(}->}[rr]\ar@{=}[dl]&&X\ar@{->>}[rr]\ar@{=}[dl]&&A\ar@{=}[dl]\\
 		\mathbb{E}:&Y\ar@{^{(}->}[rr]&&X\ar@{->>}[rr]&&A&}
 	$$
 	By hypothesis, $\mathbb{E}\varphi \in \mathcal{E}'$, and by \cite[Proposition 2.15]{Buhler}, $v'$ is an admissible monic in $\mathcal{E}$ with cokernel $K'$. By  a simple calculation, one can check that   the left-face of the upper cube is in fact a pushout diagram. Since $e$ is an admissible monic in $\mathcal{E}'$, so is  $v'$. Then, $y'=v' \circ y''$ is an admissible monic in $\mathcal{E}'$.
 	

 \end{proof}

 \subsection{Proposition.} \label{cor}
 	Let $\mathcal{E}'$ be an exact substructure of $\A$ which has enough injectives, and let $\varphi\mathcolon  B \rightarrow A$ be an $\mathcal{E}'$-phantom cover of $A$. If $\varphi$ has a kernel, then  it is $\mathcal{E}'$-injective.

 \begin{proof}
 Let  	$k: K \rightarrow B$ be a kernel of $\varphi$. Let $\mathbb{E}:\ \xymatrix{K \ar@{^{(}->}[r]^e & E \ar@{->>}[r]^a & K'}$ be a short exact sequence in $\mathcal{E}'$ with $E$ $\mathcal{E}'$-injective. Taking pushout of the short exact sequence $\mathbb{E}$ along $k$, we obtain the following diagram
 	\begin{equation}
 	\xymatrix{K \ar[r]^k\ar@{^{(}->}[d]_e&B \ar[r]^{\varphi}\ar@{^{(}->}[d]_{v}&A\ar@{=}[d]\\
 	E\ar[r]&B'\ar[r]^{\varphi'}&A.
 	}
  \end{equation}
 	By Proposition \ref{pushoutp}, $\varphi'$ is $\mathcal{E}'$-phantom. Using classical covering arguments, we deduce that  there exists a morphism $b:B' \rightarrow B$ such that $b \circ v=\id$.  From the universal property of kernel, and $e$ is an admissible monomorphism, one can easily show that $K$ is a direct summand of $E$.
 \end{proof}

 \subsection{Corollary.}\label{spec} Let $\A$ be a locally presentable category. If $\mathcal{E}'$ is  an exact substructure which is closed under direct limits and has enough injectives,   then every object of $\A$ has a $\Phi(\mathcal{E}')$-cover with $\mathcal{E}'$-injective kernel. In addition, if $\A$ has enough $\mathcal{E}'$-phantom morphisms, then  $\Phi(\mathcal{E}')$ is a special covering ideal.


  \subsection{Example.} Suppose that $\A$ is a Grothendieck category with a closed symmetric monoidal structure $\otimes$. Let $\mathcal{S}$ be a set of objects in $\A$. Then the \textit{exact structure flatly generated by} $\mathcal{S}$, denoted by $\tau^{-1}(\mathcal{S})$,  consists   of all short exact sequences $\mathbb{E}$ which remain exact under $-\otimes S$, for every $S\in \mathcal{S}$. It is clearly closed under direct limits. Besides, one can easily show that it has enough injectives.   In particular, when $\A=R\Mod$ and $S$ is a set of finitely presented $R$-modules, the exact structure  projectively generated by   $\mathcal{S}$  coincides with $\tau^{-1}(Tr(\mathcal{S}))$, where $Tr$ is  the
  Auslander-Bridger transpose of the finitely presented $R$-modules, see \cite[Theorem 8.3]{Sklyarenko}.

\subsection{Remark.} If $\A$ is locally finitely presented category,  then from $\cite{Her03}$, we know that the pure exact structure $\mathcal{E}_{\aleph_0}$ has enough injectives. If, in addition, $\mathcal{E}_{\aleph_0} \subseteq \mathcal{E}$,  then  the  ideal $\Phi(\mathcal{E}_{\aleph_0})$ is a covering ideal with pure-injective kernel. However, from Proposition \ref{proj.mor}, one can observe that if $\A$ is an abelian category without enough projective morphisms, then neither would have  enough phantom morphisms (= $\mathcal{E}_{\aleph_0}$-phantom morphisms). Then, one can not guarantee that if $\Phi(\mathcal{E}_{\aleph_0})$ is a special  covering ideal. In fact, in Section 4, we provide an example of a locally finitely presentable abelian category  in which phantom cover is just zero, in other words,  $\Phi(\mathcal{E}_{\aleph_0})=0$.

 \section{Phantom maps in $\Qco(X)$}\label{section.phantom.qc}
 In this section, we aim at investigating certain `phantom' morphisms in the category $\Qco(X)$ of quasi-coherent sheaves over a scheme $X$  which arise naturally.

 \subsection{}Given any scheme  $X$, the category $\Qco(X)$ is known to be  always a closed symmetric monoidal   Grothendieck category (for being Grothendieck, see \cite{EEr}).  So this fact  implies that it is locally $\lambda$-presentable for some regular cardinal $\lambda$, and there exist two  `pure'  exact structures:  the exact structure $\mathcal{E}_{\lambda}$  of $\lambda$-pure short exact sequences  and the exact structure $\mathcal{E}_{\otimes}$ of geometrical pure short exact sequences in $\Qco(X)$. As an aside, in most practical cases ($X$ quasi-compact and quasi-separated) the category $\Qco(X)$ is also locally finitely presented (see \cite[I.6.9.12]{GD} or \cite{Garkusha} for a precise statement).

 However, the purity that naturally arise in $\Qco(X)$ in various geometric contexts is   the so-called {\it stalk-wise} purity. We recall the following.

 \subsection{Proposition}\label{pure2} \cite[Proposition 3.4]{EEO}.
 	Let $X$ be a scheme and $\rF, \rG \in   \Qco(X)$. The following statements are equivalent:
 	\begin{itemize}
 		\item[i)] $0\ra \rF\stackrel{\tau}{\ra} \rG$ is geometrical pure exact in $\O_X \Mod$.
 		\item[ii)] There exists an open covering of $X$ by affine open sets, $\mathcal{U}=\{U_i\}$, such that $0\ra \rF(U_i)\stackrel{\tau_{U_i}}{\longrightarrow} \rG(U_i)$ is pure in $\O_X(U_i)\Mod$.
 		\item[iii)]  $0\ra \rF_x\stackrel{\tau_x}{\ra} \rG_x$ is pure in $\O_{X,x}\Mod$, for each $x\in X$.
 		
 	\end{itemize}

  \subsection{Definition.}\label{stalkwise} A short exact sequence $\mathbb E$ in $\Qco(X)$ is said to be  \textit{stalk-wise pure} if  it satisfies one of the equivalent conditions given in Proposition \ref{pure2}. We let $\mathcal{E}_{st}$ denote  the corresponding exact structure.

  \subsection{}Since stalks, colimits and tensor products in $\Qco(X)$ commute, we have the following ordering of exact structures in $\Qco(X)$
 $$\mathcal{E}_{\lambda} \subseteq \mathcal{E}_{st} \subseteq \mathcal{E}_{\otimes}.$$
As a result, there are inclusions  of relative phantom ideals $\Phi(\mathcal{E}_{\lambda}) \subseteq \Phi(\mathcal{E}_{st}) \subseteq \Phi(\mathcal{E}_{\otimes}).$ In order to understand if any of these ideals is determined by  Zariski-local property,  we make  the following definition.

 \subsection{Definition.}\label{Zariski-local.phantom}
 We call a morphism $f\mathcolon  \rG \rightarrow \rF$ in $\Qco(X)$  \emph{locally phantom} if for every  affine open subset $U$ of $X$, the induced morphism $f_U$ on  modules of sections between  $\rG$ and $\rF$ is phantom. The locally phantom morphisms in $\Qco(X)$ forms an ideal which will be denoted by  $\Phi'$.

 \vspace{3mm}

 Let us begin by proving the Zariski-local property of locally phantom morphisms.

 \subsection{Lemma.}\label{l1}
 Let $R$ be a commutative ring. A morphism $f\!:\! M'\rightarrow M$ of  $R$-modules is phantom if and only if  the induced localization morphism $f_P$ is phantom in $R_P \Mod$ for every $P \in \Spec(R)$.

 \begin{proof}
 	Suppose that $f$ is a phantom morphism in $R\Mod$. For a given prime ideal $P$ of $R$,   consider the following  pullback diagram in $R_P \Mod$
 	$$\xymatrix{0 \ar[r]& A \ar[r]\ar@{=}[d]& B' \ar[r]\ar[d]& M'_{P}\ar[d]^{f_P}\ar[r] &0\\
 		0\ar[r]&A \ar[r]& B\ar[r]&M_P \ar[r] & 0.
 	}
 	$$
 	Taking pullback over  the canonical $R$-module  morphisms $M \rightarrow M_P$ and $M' \rightarrow M'_P$, we get the following  commutative diagram in $R\Mod$
 	$$\xymatrix@!0{
 		0\ar[rr] &&A\ar[rr]\ar@{=}[dd]\ar@{=}[dr]&&N'\ar[dd]\ar[rr]\ar[dr]&&M'\ar[dd]\ar[rr]\ar[dr]^f&&0\\
 		& 0\ar[rr]&&A\ar@{=}[dd]\ar[rr]&&N\ar[dd]\ar[rr]&&M\ar[rr]\ar[dd]&&0\\
 		0\ar[rr]&&A \ar[rr]\ar@{=}[dr]&&B'\ar[rr]\ar[dr]&& M'_P \ar[rr]\ar[dr]^{f_P}&&0\\
 		&0\ar[rr]&&A\ar[rr]&&B\ar[rr]&&M_P\ar[rr]&&0
 	}$$
 	By  assumption, the upper row is  pure-exact. But purity is preserved under localization and $A_P\simeq A$ as $R_P$-modules (here $A_P$ is the localization of $R$-module $A$), see \cite[Remark 3.8]{Pinzon}. Therefore, $f_P$ is phantom in $R_P\Mod$.
 	
 	The converse follows from the fact that the localization functor preserves pullback diagrams of epimorphisms and  that  a monomorphism $\iota$ is pure if and only if  for every $P \in \Spec(R)$,  the induced localization morphism  $\iota_P$  is pure.
 \end{proof}

 \subsection{Theorem.}\label{Zariski-local.phantom.theorem}
 	Let $f\mathcolon \rG \rightarrow \rF$ be a morphism in $\Qco(X)$. The following are equivalent:
 	\begin{itemize}
 		\item[i)] $f$ is locally phantom.
 		\item[ii)] There is a cover $\mathcal{U}$ of $X$ consisting of affine open subsets such that $f_U$ is phantom for every $U \in \mathcal{U}$.
 		\item[iii)] $f_x$ is phantom in $\O_{X,x}\Mod$, for all $x \in X$.
 	\end{itemize}

 \begin{proof}
 	Follows from Lemma \ref{l1}.
 \end{proof}

 \subsection{Proposition.}\label{ph1}
 Let $f\mathcolon \rG \rightarrow \rF$ be a morphism in $\Qco(X)$. If $f$ is locally phantom  then it is also $\mathcal{E}_{st}$-phantom. That is, $\Phi'\subseteq \Phi(\mathcal{E}_{st})$.

 \begin{proof}
 	Let
 	$$\xymatrix{\mathbb{E}f: \textrm{ }0\ar[r] & \rF'\ar[r]\ar@{=}[d] & \rG'\ar[r]\ar[d] & \rG\ar[r]\ar[d]^f & 0\\
 		\mathbb{E}: \textrm{ }0\ar[r] & \rF'\ar[r] & \rF''\ar[r] & \rF\ar[r] & 0
 	}$$
 	be a pullback diagram in $\Qco(X)$. For every  affine open subset $U \subseteq X$,  the sequences $\mathbb{E}(U)$ and  $(\mathbb{E}f)(U)$ of modules of sections over $U$ are  exact and  they  induce a pullback diagram in $\O_X(U) \Mod$.
 	By assumption, $(\mathbb{E}f)(U)$ is pure-exact, and therefore, $\mathbb{E}\in \mathcal{E}_{st}$.
 \end{proof}

\subsection{Proposition.}\label{ph2}
	If  $X$ is semi-separated, then $f\mathcolon  \rG \rightarrow \rF$ is  $\mathcal{E}_{st}$-phantom in $\Qco(X)$ if and only if it is locally phantom.  Moreover, in this case, $\Phi(\mathcal{E}_{\otimes})=\Phi(\mathcal{E}_{st})=\Phi'$.

\begin{proof}
	For  an affine open subset $U$, let
	$$\xymatrix{0\ar[r] & N\ar[r]\ar@{=}[d] & M'\ar[r]\ar[d] & \rG(U)\ar[r]\ar[d]^{f_U} & 0\\
		0\ar[r] & N\ar[r] & M\ar[r] & \rF(U)\ar[r] & 0
	}$$
	be a pullback diagram in $\O_X(U)\Mod$. Note that $\iota_* $ is an exact functor from $\Qco(U)$ to $\Qco(X)$ since $X$ is semi-separated. Consider the following commutative diagram
	$$\xymatrix@!0{
		0\ar[rr] &&\iota_{\ast}\widetilde{N}\ar[rr]\ar@{=}[dd]\ar@{=}[dr]&&\rT'\ar[dd]\ar[rr]\ar[dr]&&\rG\ar[dd]\ar[rr]\ar[dr]^f&&0\\
		& 0\ar[rr]&&\iota_{\ast}\widetilde{N}\ar@{=}[dd]\ar[rr]&&\rT\ar[dd]\ar[rr]&&\rF\ar[rr]\ar[dd]&&0\\
		0\ar[rr]&&\iota_{\ast}\widetilde{N}\ar[rr]\ar@{=}[dr]&&\iota_{\ast}\widetilde{M'}\ar[rr]\ar[dr]&&\iota_{\ast}\rG|_U\ar[rr]\ar[dr]&&0\\
		&0\ar[rr]&&\iota_{\ast}\widetilde{N}\ar[rr]&&\iota_{\ast}\widetilde{M}\ar[rr]&&\iota_{\ast}\rF|_U\ar[rr]&&0
	}$$
	where each face is a pullback diagram in $\Qco(X)$. By assumption, the upper exact sequence is stalkwise pure. Then the short exact sequence of modules of  sections over $U$,  $0 \rightarrow N \rightarrow \rT'(U) \rightarrow \rG(U) \rightarrow 0$, which is   isomorphic to $0 \rightarrow N \rightarrow M' \rightarrow \rG(U)  \rightarrow 0$,
	is pure-exact. So $f_U$ is a phantom morphism. This shows that $\Phi'=\Phi(\mathcal{E}_{st})$. Finally, by \cite[Proposition 2.10]{EGO}, we have that $\mathcal{E}_{\otimes}=\mathcal{E}_{st}$ for any quasi-separated scheme (so, in particular, for any semi-separated scheme). Therefore, $\Phi'=\Phi(\mathcal{E}_{st})=\Phi(\mathcal{E}_{\otimes})$.
\end{proof}

 \subsection{Lemma.}
 	The ideals $\Phi(\mathcal{E}_{st})$, $\Phi(\mathcal{E}_{\otimes})$ and $\Phi'$ in $\Qco(X)$ are closed under direct limits. The ideal $\Phi(\mathcal{E}_{\lambda})$ is closed under direct limits when $\lambda=\aleph_0$, that is, when $\Qco(X)$ is locally finitely presented.

 \begin{proof}
 	It follows from the fact that $\mathcal{E}_{st}$ and $\mathcal{E}_{\otimes}$ are closed under direct limits and finite limits commute with direct limits. Finally, note that the result  holds for $\mathcal{E}_{\lambda}$ when $\lambda=\aleph_0$.
 \end{proof}

 \subsection{Corollary.}\label{covering.ideals.qco}
 	The following holds in $\Qco(X)$:
 	\begin{enumerate}
 		\item The  ideals $\Phi(\mathcal{E}_{st})$ and $\Phi(\mathcal{E}_{\otimes})$  are covering ideals  with stalkwise pure-injective and geometrical pure injective kernel, respectively.
 		\item The ideal $\Phi'$ is covering.
 	\end{enumerate}
In addition, if $X$ is a quasi-compact and semi-separated scheme, then $\Phi(\mathcal{E}_{st})= \Phi(\mathcal{E}_{\otimes})=\Phi'$ is a special covering ideal.
\begin{proof}
	From \cite[Theorem 4.10]{EEO} and \cite{EGO15}, we know that $\mathcal{E}_{st}$ and $\mathcal{E}_{\otimes}$ have enough injectives.
	
	One can easily check that any morphism $f\mathcolon \rF\to \rG$ in $\Qco(X)$, with $\rF$ a flat quasi-coherent sheaf, belongs to  $\Phi'$ (and thus, to $\Phi(\mathcal{E}_{st})\subseteq \Phi(\mathcal{E}_{\otimes})$). If $X$ is a quasi-compact and semi-separated scheme,  $\Qco(X)$ has a flat generator,   see  \cite{Leo}. Then apply Corollary \ref{spec}.
	
\end{proof}

 Now, we claim that there are no non-zero $ \mathcal{E}_{\lambda}$-phantom morphisms in $\Qco(X)$  even over nice non-afine schemes. For that, we firstly prove the following.

 \subsection{Proposition.}\label{proj.mor}
 	
 	Let $(\A,\mathcal{E})$ be an exact category. Assume $\A$ to be  locally finitely presented and $\mathcal{E}_{\aleph_0} \subseteq \mathcal{E}$. Given a morphism $\phi: M \rightarrow N  \in \A$, $\phi\in \Phi(\mathcal{E}_{\aleph_0})$ if and only if it is a direct limit of projective morphisms.

 \begin{proof} Suppose that $\phi\in \Phi(\mathcal{E}_{\aleph_0})$. We may write $M=\varinjlim M_i$ as a direct limit of finitely presented objects. Let $\{\tau_i\mathcolon   M_i\to M\}$ be the structural morphisms. It is easy to check that the family $\{\phi\circ \tau_i\}$ is a directed system of projective morphisms and that $\phi=\varinjlim(\phi\circ\tau_i)$.

 	Conversely assume that $\{\phi_i \mathcolon M_i\to N_i\}_{i \in I}$ is a  directed system of  a projective morphisms with  $\phi=\varinjlim \phi_i$. Let $\mathbb{E}$ be a short exact sequence in $\mathcal{E}$ ending with $N$. Consider the following pullback diagram in $\A$
 	$$\xymatrix{\mathbb{E}\phi: \textrm{ } & K \ar@{=}[d]\ar@{^{(}->}[r] & X'\ar@{->>}[r]\ar[d] &\varinjlim M_i \ar@{->}[d]^{\varinjlim \phi_i} \\
 		\mathbb{E}: \textrm{ }  & K\ar@{^{(}->}[r] & X \ar@{->>}[r] &  \varinjlim N_i.
 	}$$
 Let $T\to \varinjlim M_i$ be a morphism with $T$ finitely presented. Then $T\to\varinjlim M_i $ factors through $T\to M_i\to \varinjlim M_i$ for some  $i\in I$. Taking pullback of $\mathbb{E}$ through the structural morphism $\kappa_i\mathcolon N_i\to \varinjlim N_i$, we have
 	$$\xymatrix{   K \ar@{=}[d]\ar@{^{(}->}[r] & X_i'\ar@{->>}[r]\ar[d] &N_i \ar@{->}[d]^{\kappa_i} \\
 		  K\ar@{^{(}->}[r] & X \ar@{->>}[r] &  \varinjlim N_i.
 	}$$
 As $\phi_i\mathcolon M_i\to N_i$ is projective, there exists a  morphism $M_i\to X_i'$ making the obvious diagram commutative. Finally, since the $\mathbb{E}\phi$ is a pullback, there is a morphism $T\to X'$ such that $T\to X'\to \varinjlim M_i$ equals to $T\to \varinjlim M_i$. Hence, the short exact sequence $\mathbb{E}\phi$ is $\aleph_0$-pure.
 	
 \end{proof}

 \subsection{} Let $X={\bf P^1}(R)$ be  the projective line over any commutative ring $R$.  Let us cover $X$ by the usual affine open subsets $\mathcal U= U\hookleftarrow U\cap V\hookrightarrow V$. The structure sheaf of $X$ is given by the following representation of $\mathcal U$,
 $$\mathcal O= R[x]\hookrightarrow R[x,x^{-1}]\hookleftarrow R[x^{-1}].$$

 The Serre's twisting sheaves $\mathcal O(n)$ are given by
 $$\mathcal O(n)=R[x]\hookrightarrow R[x,x^{-1}]\stackrel{x^n}{\fiz}R[x^{-1}],$$
 with
 $n\in \Z$. It is known that the family of twisting sheaves $\{\O(n)\}_{n\in \Z}$ generates the category $\QcoP$. Indeed it suffices to take the family $\{\O (-n)\}_{n\in\N}$ to generate $\QcoP$.

 On the other hand, any quasi-coherent sheaf $\M\in \Qco({\bf P^1}(R))$ is determined by a representation of $\mathcal U$
 $$\T:\quad M\stackrel{f}{\to} P\stackrel{g}{\fiz} N,$$
 where  $M\in \Re\Mod$,
 $N\in \Ri\Mod$, $P\in \Ro\Mod$,  $f$ is an $\Re$-linear map
 and $g$, a $\Ri$-linear such that
 $$\Si f\!:\! \Si M\to  P \quad \textrm{ and } \quad  \B g\!:\!\B N\to P$$
 are isomorphisms, where $S=\{1,x,x^2,\cdots\}$ and
 $T=\{1,x^{-1},x^{-2},\cdots\}$.

 \subsection{Lemma.}\label{prev}
 	Let $0\neq \T=(M\stackrel{f}{\to}P\stackrel{g}{\fiz}N)\in \Qco({\bf P^1}(R))$.
 	Given $0\neq m\in M$ (resp., $0\neq p\in P$, $0\neq y\in N$), there exists a
 	natural number
 	$k_m$ (resp., $k_p$, $k_y$) such that for every $l\geq k_m$ (resp., $l\geq k_p$, $l\geq k_y$) and every subset $\Delta\subseteq \mathbb Z$, any morphism
 	$$(\gamma_1,\gamma,\gamma_2)\mathcolon  \T \to \oplus_{n\in \Delta}\O(-n-l)$$
 	maps $m$ (resp., $p$, $y$) to
 	zero.

 \begin{proof}  Since $T^{-1}g\mathcolon T^{-1}N\to P$ is an
 	isomorphism, we get that
 	$f(m)=g(a)/x^{-l}$, for some $a\in N$ and $l\in \N$. That is, $g(a)=x^{-l}f(m)$. Set $k=l+1$.
 	
 	Let $\Delta\subseteq \mathbb Z$ and consider a morphism $(\gamma_1,
 	\gamma,\gamma_2):\T\to \oplus_{n\in \Delta}\O(-n-k)$. Let us write
 	$$\gamma_1(m)=(\ldots,p_1(x),\cdots,p_k(x),\ldots)$$ and
 	$$\gamma_2(a)=(\ldots,q_1(x^{-1}),\cdots,q_t(x^{-1}),\ldots).$$
 	Then,
 	$$\gamma \circ g(a)=\gamma(x^{-l}f(m))=x^{-l}\gamma(f(m))=
 	\gamma_1(m)=x^{-l}(p_1(x),\cdots,p_k(x)).$$ Thus,
 	$ord(x^{-l}p_i(x))\geq -l$, for every $1\leq i\leq k$. But, by the
 	commutativity of the diagram, we also get that
 	$$\gamma\circ g(a)=\oplus_{n\in \Delta} x^{-(n+k)}(\gamma_2(a))=\oplus_{n\in \Delta}
 	x^{-(n+k)}(q_1(x^{-1}),\cdots,q_t(x^{-1}))$$ $$
 	=x^{-k}(r_1(x^{-1}),\cdots,r_t(x^{-1})),$$ with
 	$ord(x^{-k}r_i(x^{-1}))\leq -k=-l-1$, for all $1\leq i\leq t$. This shows that
 	$\gamma_1(m)=0$.
 \end{proof}

 \subsection{Corollary.}\label{nomorfp}
 	Assume that $\T\in \Qco({\bf P^ 1}(R))$ is finitely presented. Then, there exists a
 	natural number $k=k(\T)$ such that there are no nonzero morphisms from $\T$
 	into an arbitrary direct sum of the elements of the family $\{ \O(-n-k):\ n\in\N
 	\}$.

 \begin{proof} Without lost of generality, we may assume that $M$ is generated by $m_1,\ldots, m_s$;
 	$P$ is generated by $p_1,\ldots, p_s$; and $N$ is generated by $y_1,\ldots,y_s$. Then,
 	in view of Lemma \ref{prev}, we just have to take $k\geq
 	\textrm{max}\{k_{m_i},k_{p_i},k_{y_i}:\ i=1,\ldots, s\}$.
 \end{proof}
 We can argue now that there are no nonzero projective morphisms in $\QcoP$. The argument is based on the same strategy followed in  \cite[Corollary 2.3]{EEGO} to show that there are no nonzero projective objects in $\QcoP$.
 \subsection{Proposition.}\label{nohay}
 	The only projective morphism in $\Qco({\bf P^1}(R))$ is the zero map.

 \begin{proof}
 	Suppose that $\M\to\L$ is a projective morphism in $\QcoP$. Since $\QcoP$ is locally finitely presented, we can assume that $\M$ is finitely presented. There exists an epimorphism $\oplus_{n\in
 		\N} \O(-n)^{(X_n)}\to \L$. Let $(r,s,v)_{-n}$ be the composition map
 	$$\O(-n)\hookrightarrow \oplus_{n\in
 		\N} \O(-n)^{(X_n)}\to \L.$$
 	
 	Let us fix a natural number $n_0$. For any $n\in \mathbb N$, we can consider the morphisms $\O(-n-n_0)\to
 	\L$ given by $(r,s,x^{-n_0}v)_{-n-n_0}$ and $(x^{n_0}r,x^{n_0}s,v)_{n+n_0}$, respectively. These
 	two morphisms induce a morphism
 	$$\O(-n-n_0)\oplus \O(-n-n_0)\to \L.$$
 	
 	In turn, these morphisms induce an epimorphism
 	$$\label{eqon}\bigoplus_{n\in \N}\left (
 	\O(-n-n_0)^{(X_n)}\oplus
 	\O(-n-n_0)^{(X_n)}\right )\stackrel{\psi}{\longrightarrow} \L.
 	$$
 	Now, as $\M\in\QcoP$ is finitely presented,
 	Corollary
 	\ref{nomorfp} states that there exists an $n_0\in \N$ such that it is not
 	possible to factorize any nonzero morphism $\M\to \L$ through ${\psi}$.
 \end{proof}

 \subsection{Theorem.}\label{nohay.phantom}
 	There are no non-trivial phantom maps in $\QcoP$.

 \begin{proof}
 	The category $\QcoP$ is locally finitely presented, so the result follows from Propositions \ref{nohay} and \ref{proj.mor}.
 \end{proof}

\end{document}